\documentclass[10pt]{amsart}
\usepackage{amsmath,amssymb,bbm}
\usepackage{amsthm}
\usepackage[all]{xy}

\newcommand{\CC}{{\mathbb C }}
\newcommand{\PP}{{\mathbb P }}
\newcommand{\QQ}{{\mathbb Q }}
\newcommand{\RR}{{\mathbb R }}
\newcommand{\ZZ}{{\mathbb Z }}

\newcommand{\II}{{\mathsf I }}
\newcommand{\JJ}{{\mathsf J }}

\newcommand{\TT}{{\mathsf T }}

\newcommand{\cA}{\mathcal{A}}

\newcommand{\dz}{{\partial_z }}
\newcommand{\dt}{{\partial_t }}
\newcommand{\dw}{{\partial_w }}
\newcommand{\dx}{{D}}

\newcommand{\us}{{u^*}}
\newcommand{\bu}{{\mathbf u}}

\newcommand{\gl}{\lambda}
\newcommand{\gL}{\Lambda}

\newcommand{\detr}{\mathrm{det}_{\mathrm{right}}}
\newcommand{\Spec}{\mathrm{Spec}\:}
\newcommand{\res}{\mathrm{res}}

\newcommand{\prf}{\rm\noindent\textbf{Proof. }}
\newcommand{\qd}{\hfill$\blacksquare$}

\newcommand{\lqs}{\leqslant}
\newcommand{\gqs}{\geqslant}

\newcommand{\IS}{$1^{\mathrm{st}}$}
\newcommand{\RS}{$2^{\mathrm{nd}}$}

\newcounter{pphcounter}[section]
\renewcommand{\thepphcounter}{\thesection.\arabic{pphcounter}}
\newcommand{\pph}[1]{\noindent
\refstepcounter{pphcounter}\bf\thepphcounter. #1\rm}

\begin{document}
\bibliographystyle{alpha}

\title{Fuchsian equations of type DN}

\author{Vasily Golyshev and  Jan Stienstra}
\address{Vasily Golyshev,
Number Theory Section, Steklov Mathematical Institute,
Russian Academy
of Sciences, 8 Gubkina str., Moscow 119991, Russia}

\address{Jan Stienstra, Department of Mathematics, Utrecht University,
Postbus 80.010,$\qquad$
3508 TA, Utrecht, The Netherlands}

%%\date{\today}

\maketitle

\parbox{340pt}{\small \bf Abstract. \rm We prove that a generic differential
operator of type DN is irreducible, regular, (anti)self-adjoint, and
has quasiunipotent local monodromies.
We prove that
the defining matrix of a DN operator can be recovered from the
expression of the operator as a polynomial in $t$ and $\dt$.}

\bigskip
\bigskip
\bigskip

\section*{Introduction.}
\label{introduction}
Let $A=(a_{ij})_{0\lqs i,j\lqs N}$ be a matrix with entries in $\CC$
satisfying
$a_{ij}=0$ if $i-j>1$ and $a_{ij}=1$ if $i-j=1$ and
$a_{ij}=a_{N-j,N-i}$.
With this matrix we associate a differential operator $L_{A,\infty}$
as follows. Let
$\tilde{A}=(a_{ij} \dt^{j-i+1})$ and $\dt=\frac{d}{dt}$. The right determinant
$
\detr (t\dt -\tilde{A})
$
is a differential operator which is uniquely divisible from the right by $\dt$. We define
$$
L_{A,\infty}\,=\,\detr ( t\dt  -\tilde{A} )\dt^{-1} \,.
$$
The notation $L_{A,\infty}$ signifies its dependence on $A$ and the
fact that the space of its solutions has maximally unipotent
monodromy at infinity. For $\Phi$ in a given module over the ring of
the differential operators (e. g. that of functions of the variable
$t$), the differential equation $L_{A,\infty} \Phi =0$ was called
\emph{the determinantal equation of order $N$}, \rm or simply
equation DN, in \cite{Go}. Theorem  \ref{thm:DN--} in the present
paper characterizes (through a constructive bijective
correspondence) these operators as precisely those differential
operators of the form
$$
(t\dt)^{N}t + \sum_{p = 1}^{N + 1} g_p(t\dt )\,\dt^{p-1}\,
$$
with $g_p$ a polynomial of degree $\lqs N-p+1$ in $t\dt$, such that
$
g_p(\mathrm{argument})=(-1)^{N-p+1}g_p(-\mathrm{argument}-p).
$

The construction of the differential operators of type DN in
\cite{Go} is motivated by  mirror symmetry for minimal Fano
varieties (i.e.  ones whose cohomology is just $\ZZ$
in every even dimension), or complete intersections therein.
One expects to obtain a
Picard--Fuchs equation in the Landau--Ginzburg model of a given Fano
when one specializes $a_{ij}$ to its two--pointed Gromov--Witten
invariants in a standard way \cite{Go}. We remark here that this
construction was essentially introduced by Dubrovin under the name
of `second structural connection'. For a textbook description of
Dubrovin's theory we refer the reader to \cite{Ma}.

There is, however, no general  characterization of the
differential operators of type DN (or matrices $A$) that actually arise from the
enumerative geometry of Fano varieties.
A natural question, in the light of mirror symmetry, is then
to ask which of the above differential operators are Picard-Fuchs operators
for 1-parameter families of complex varieties, or more generally variations of
 Hodge structure.

According to Dubrovin (see \cite{Ma} Chapter II)
a generic differential operator DN is associated
with a Fuchsian connection, hence is regular. We prove that it is irreducible
and coincides,
up to a sign, with its adjoint.
It is known to have maximal unipotent monodromy at $t=\infty,$
and we show in this note that it has quasi-unipotent
monodromies elsewhere. Regularity and quasi-unipotence of local monodromies
are definite requirements for Picard-Fuchs operators by Deligne,
but do not suffice.

Unfortunately, we know no general way to tell  which of the differential
operators of type DN do also satisfy further requirements for coming
from a variation of $\RR$-Hodge structures or being globally crystalline.

\

\noindent \textbf{Acknowlegdment.} The authors thank NWO and RFBR
(NWO-RFBR grant  047.017.019, NWO-RFBR grant 047.011.2004.026 [RFBR
05-02-89000-NWO\_\,a] for Dutch Russian Research Cooperation) and
the Mathematics Department of Utrecht University for financial
support for visits, during which most research for this paper was
done. V.Golyshev was supported by \mbox{grant  RUM1-2661-MO-05 } of
the U.S. Civilian Research and Development Foundation for the
Independent States of the Former Soviet Union (CRDF).

\section{Algebraic preliminaries.}
\label{section:general algebra}

\pph{Conventions.}\label{conventions}
In this section $R$ is an associative ring with unity.
We do not want to use ugly looking indices such as $a_{iN}.$
For this reason we set $n=N,$ and produce our DN's starting
with a matrix of size $n+1$
whose row and column indices run from $0$ to $n$.

\

\pph{Definition.} \it \label{def:detright}
The \emph{right determinant} $\detr (M)$ of a square matrix
$M=(M_{ij})_{0\lqs i,j\lqs n}$
with entries in $R$ is defined by
expanding with respect to the right column:
$$
\detr(M)=\sum_{i=0}^n M_{in} C_{in}
$$
where
$C_{in}$ is the cofactor the element
$M_{in}$; the cofactor $C_{in}$ is in turn a right
determinant times a sign $(-1)^{i+n}$.

\rm
By fully expanding this recursive definition one sees that
$$
\detr (M)=\sum_\sigma\mathrm{sign}(\sigma)M_{\sigma(n),n}\cdots
M_{\sigma(0),0}
$$
where $\sigma$ runs over all permutations of $\{0,\ldots, n\}$.
In particular, if $M$ has a row of $0$'s, then $\detr (M)=0$.

\

\pph{Proposition.} \label{diag-change-signs} \it If two matrices
$M=(M_{ij})_{0\lqs i,j\lqs n}$ and $M'=(M'_{ij})_{0\lqs i,j\lqs n}$ are related by  $M'_{ij}=(-1)^{j-i+1}M_{ij}$,
then $\detr(M')=(-1)^{n+1}\detr (M)$. \rm

\qd

\

\pph{Definition.} \it \label{def:almost triangular}
Let $M=(M_{ij})$ be a square matrix with entries
in $R$. We say that $M$ is \emph{almost triangular}
  \footnote{{In numerical linear algebra, a matrix $M$ is said to be
   ``upper (resp. lower) almost triangular"
   or in ``upper (resp. lower) Hessenberg form''   if $M_{ij}=0$
   for $i > j+1$ (resp.  $j > i+1$).  Our almost triangular matrices
are thus the ``upper almost triangular" matrices of numerical linear
algebra subject to the additional requirement that
all the subdiagonal elements be $-1$'s.}}
if $M_{ij} = 0$ for $i > j + 1$, and $M_{j + 1, j} = - 1$.

\

\pph{}\label{detright1}
For an almost triangular matrix $M$ one can reformulate \ref{def:detright}
as a simple inductive algorithm (cf. \cite{GR}) that consecutively
expresses a principal minor in terms
of preceding ones:
$$
 P_0 = 1, \quad P_{j + 1} = \sum^j_{i = 0} M_{ij} P_i,\;
\detr ( M ) = P_{n + 1}  .
$$
More generally, $P_0,\ldots,P_{n + 1}$ are the right principal minors of $M$.

\

\pph{Proposition.}  \label{detright2} {\it
Let $M$ be an almost triangular matrix of size $n$. Define the elements $Q_j$
in $R$ inductively by:
$$
 Q_0 = 1, \quad Q_{j + 1} = \sum^j_{i = 0}  Q_i M_{n-j,n-i}.
$$
Then $Q_{n+1}=\detr (M)$.}

\

\prf
When fully expanded both inductive algorithms lead to:
$$
P_{n + 1} = Q_{n+1} = \sum M_{(i_k,i_{k+1}-1)}M_{(i_{k-1},i_k-1)}
\cdot\ldots\cdot M_{(i_0,i_1-1)}
$$
where the sum is over all sequences
$(i_0,\ldots,i_{k+1})$ of integers satisfying
$0=i_0<i_1<\ldots<i_k<i_{k+1}=n+1$.
\qd

\

\pph{Definition.}\label{def:transpose} \it
For a square matrix $M=(M_{ij})_{0\lqs i,j\lqs n}$ with entries
in $R$
let $M^\tau$ be the matrix with $(i,j)$-entry $M^\tau_{ij}=M_{n-j,n-i}$.

\medskip

\rm Thus $M^\tau$ is the `transpose of $M$ with respect to the anti-diagonal'.
It relates to the ordinary transpose $M^t$ as
$M^\tau=\JJ M^t\JJ$, where $\JJ=(\JJ_{ij})_{0\lqs i,j\lqs n}$ denotes the
matrix with $\JJ_{ij}=1$ if $i+j=n$ and
$\JJ_{ij}=0$ otherwise.

\

\pph{Proposition.} \label{prop:CH} {\it
Let $R$ be as above and let $E$ be a right $R$-module.
Let $\xi_0, \ldots, \xi_n$ be elements of $E$.
Let $M=(M_{ij})_{0\lqs i,j\lqs n}$ be an almost triangular matrix over $R$.
Then}
$$
( \xi_0, \ldots, \xi_n )M^\tau = (0,\ldots,0)
\qquad\Longrightarrow\qquad
\xi_0\detr ( M ) = 0.
$$

\

\prf
If $( \xi_0, \ldots, \xi_n )M^\tau = (0,\ldots,0)$, then
$\xi_{j+1}=\sum_{i=0}^j \xi_i M_{n-j,n-i}$  for $j=0, \dots, n-1$.
This implies  $\xi_j=\xi_0 Q_j$ for
$0\lqs j\lqs n$ with $Q_j$ as in \ref{detright2}.
Now $( \xi_0, \ldots, \xi_n )(\text{right column of } M^\tau) = 0$
implies
$\xi_0\detr ( M ) = 0$.
\qd

\

\pph{}
Suppose that the ring $R$ is equipped an anti-involution ${}^\vee$;
i.e. a map $R\rightarrow R\,,\; x\mapsto x^\vee$, such that
$(x+y)^\vee=x^\vee+y^\vee$ and $(xy)^\vee=y^\vee x^\vee$ for all $x,y\in R$.
Let $M=(M_{ij})$ be an almost triangular matrix over $R$.
Applying the anti-involution to the inductive algorithm in \ref{detright2}
we obtain
$$
 Q_0^\vee = 1, \quad Q_{j + 1}^\vee = \sum^j_{i = 0}  M_{n-j,n-i}^\vee Q_i^\vee
\,,
$$
which is in fact the inductive algorithm in \ref{detright1} for computing
the right determinant of the matrix $M^{\tau\vee}$
(on matrices ${}^\vee$ acts componentwise).
Thus

\

\pph{Proposition.}\label{anti-involution} \it One has
$
(\detr (M))^\vee\,=\,\detr(M^{\tau\vee})\,.
$ \rm

\qd

\section{Almost triangular matrices over the Weyl algebra.}
\label{section:special algebra}

\pph{}\label{operator expansion}
Let $B=\QQ[a_{ij}]$ be the
commutative polynomial ring in the variables $a_{ij}$ with
$0\lqs i\lqs j\lqs n$.
We put a grading on $B$
such that $a_{ij}$ is homogeneous of degree $j-i+1$.
Let $R$ be the \emph{Weyl algebra over $B$}, i.e.
the non-commutative polynomial ring $B[u,\us]$
with centre $B$ modulo the commutation relation $u\us-\us u=1$.
We define the matrix $\tilde{A}=(\tilde{A}_{i j})_{0\lqs i,j\lqs n}$ with entries in
$R$ by
$\tilde{A}_{i j}=0$ if $i>j+1$, $\tilde{A}_{i j}=1$ if $i=j+1$ and
$$
\tilde{A}_{ij} = a_{i j} u^{j - i + 1}\quad \textrm{if}\quad i<j+1\,.
$$
Then $u\us-\tilde{A}$ is an almost triangular matrix; here, and henceforth,
we simplify the notation by writing just $u\us$ instead of $u\us$
times the identity matrix of size $n+1$.
Using the inductive algorithm in \ref{detright1} one checks that
its right determinant has an expansion of the form:

\

\pph{}\label{expansion}
$\qquad\displaystyle
\detr ( u\us -\tilde{A} ) =
(u\us)^{n + 1} + \sum_{p = 1}^{n + 1} \sum_{k=0}^{n-p+1}
x^{(p)}_k u^p (u\us)^k
$\\
in which $x^{(p)}_k$ is a homogeneous element of degree $p$ in $B$.

\

\pph{}\label{tau}
Define, for $p\gqs 1$, the endomorphisms $\tau^{\gqs p}$ and $\tau^{\lqs p}$ of the ring $B$ by
$$
\begin{array}{llllll}
\tau^{\gqs p}(a_{ij})=0 & \textrm{if} &  j-i+1<p\,,&
\tau^{\gqs p}(a_{ij})=a_{ij} & \textrm{if} &  j-i+1\gqs p\,,
\\
\tau^{\lqs p}(a_{ij})=0 & \textrm{if} &  j-i+1>p\,,&
\tau^{\lqs p}(a_{ij})=a_{ij} & \textrm{if} &  j-i+1\lqs p\,.
\end{array}
$$
Then $\tau^{\lqs p}(x^{(p)}_k)=x^{(p)}_k$
and $\tau^{\gqs p}(x^{(p)}_k)$ is a linear combination of
the variables $a_{i,p+i-1}$, since $x^{(p)}_k$ is homogeneous of degree $p$.
Thus for every $p$, $1\lqs p\lqs n$, there is a square
matrix $K^{(p)}=(K^{(p)}_{ki})$ over $\QQ$ of size $n+2-p$ such that
$$
\tau^{\gqs p}(x^{(p)}_k)=\sum_{i=0}^{n-p+1} K^{(p)}_{ki}\,a_{i,p+i-1}\,.
$$

\

\pph{Lemma.}\label{invertible} {\it
The matrix $K^{(p)}$ is invertible over $\QQ$, for every $p\gqs 1$.}

\

\prf For $j\gqs i\gqs 0$, let $E^{(i,j)}=(E^{(i,j)}_{kl})_{0\lqs k,l\lqs n}$
denote the matrix with $E^{(i,j)}_{ij}=u^{j-i+1}$, $E^{(i,j)}_{k+1,k}=1$ for
$0\lqs k\lqs n-1$ and $E^{(i,j)}_{kl}=0$ else. The recursion rule
in \ref{detright1} yields
\begin{eqnarray*}
\detr ( u\us - E^{(i,p+i-1)} ) - (u\us)^{n + 1}
&=& - (u\us)^{n+1-p-i} u^p (u\us)^i \\
=\; -u^p (u\us-p)^{n+1-p-i}  (u\us)^i
&=&
\sum_{k=0}^{n-p+1} K^{(p)}_{ki} u^p (u\us)^k\,.
\end{eqnarray*}
Invertibility of the matrix $K^{(p)}$ therefore follows from the fact that
the elements
$(u\us-p)^{n+1-p-i}  (u\us)^i$ with $0\lqs i\lqs n-p+1 $ are
linearly independent.
\qd

\

\pph{Theorem.}\label{reconstruction thm} {\it  Consider the polynomial ring $\Lambda=\QQ[\lambda^{(p)}_k]_{0\lqs k\lqs n-p+1\lqs n}$
together with the ring morphism $\varphi:\Lambda\rightarrow B\,,
\; \varphi(\lambda^{(p)}_k)=x^{(p)}_k$.
Then $\varphi$ is an isomorphism.}

\

\prf
From \ref{tau} and \ref{invertible} we see that $a_{i,p+i-1}$ is a linear combination of the elements $x^{(p)}_k$ with $0\lqs k\lqs n-p+1$
plus a polynomial in the elements $a_{ij}$ with $j-i+1<p$. So, by induction,
$a_{i,p+i-1}$ is in the image of $\varphi$. Hence $\varphi$ is surjective.
Put a grading on $\Lambda$ by declaring that $\lambda^{(p)}_k$ is homogeneous
of degree $p$. Then $\varphi$ is a morphism between graded
$\QQ$-algebras. Moreover for every $p$ the homogeneous pieces of degree
$p$ in $\Lambda$ and $B$ are $\QQ$-vector spaces of the same dimension.
So in each degree $\varphi$ is a surjective linear map between
vector spaces of the same dimension. Hence $\varphi$ is an isomorphism.
\qd

\

In down to earth terms the theorem means:

\

\pph{Corollary.}\label{reconstruction} {\it
The matrix $A$ can be reconstructed from the expansion of\\
$\detr ( u\us -\tilde{A} )$.}

\

\prf
By the theorem the coefficients $x^{(p)}_k$ in the expansion of
$\detr ( u\us -\tilde{A} )$ are algebraically independent and the matrix
entries $a_{ij}$ are polynomials in the $x^{(p)}_k$'s.
This remains true when the $x^{(p)}_k$'s
are specialized to complex numbers, leading to an almost triangular matrix
$-A$ over $\CC$.
\qd

\section{Realizations}
\label{section:Realizations}

\pph{}\label{realization problem}
In this section $R$ is the \emph{Weyl algebra over $\CC$}, i.e.
the non-commutative polynomial ring $\CC[u,\us]$
with centre $\CC$ modulo the relation $u\us-\us u=1$.
Further, $A=(a_{ij})_{0\lqs i,j\lqs n}$ is a matrix with entries in $\CC$
such that $-A$ is almost triangular.
We define the matrix $\tilde{A}=(\tilde{A}_{i j})_{0\lqs i,j\lqs n}$ with
entries in $R$ by
$$
\tilde{A}_{ij} = a_{i j} u^{j - i + 1}\,.
$$
We want to apply the result of Proposition \ref{prop:CH} to the almost triangular matrix
$u\us-\tilde{A}^\tau$. So, we need a right $R$-module $E$ and in it elements
$\xi_0, \ldots, \xi_n$ such that
$$
( \xi_{0, \ldots,} \xi_n )(u\us-\tilde{A}) = (0,\ldots,0)\,.
$$
We will present two realizations of this situation, called the \IS-model and the \RS-model.
The terminology \IS-model and \RS-model refers to the fact that these
give the \emph{first} and the \emph{second structure connection}, respectively,
in \cite{Ma} Ch. II; see \ref{IS-model(bis)} and \ref{RS-model} below.
For both realizations we use an isomorphism of the Weyl algebra with
an algebra of differential operators:
$$
\begin{array}{lllll}
\textrm{\IS-model:}&\quad R\stackrel{\sim}{\rightarrow} \CC[z,\dz]\,,&
u\mapsto z\,,& \us\mapsto -\dz&\quad (\dz=\frac{d}{dz}),
\\
\textrm{\RS-model:}&\quad R\stackrel{\sim}{\rightarrow} \CC[t,\dt]\,,&
u\mapsto \dt\,,& \us\mapsto t&\quad (\dt=\frac{d}{dt}).
\end{array}
$$
These isomorphisms with the Weyl algebra yield the isomorphism
$$
\CC[t,\dt]\stackrel{\sim}{\rightarrow} \CC[z,\dz]\,,\qquad
t\mapsto -\dz\,,\; \dt\mapsto z\,.
$$
This means that the
\IS-model as a left $\CC[z,\dz]$-module is the Fourier transform of the \RS-model
as a left $\CC[t,\dt]$-module in the standard sense
(cf. \cite{Ka} p.71).

\

\pph{}\label{leftright}
The algebra $\CC[z,\dz]$ admits an anti-involution
${}^\vee$ which is the identity on $\CC$ and  satisfies
$$
z^\vee=z\,,\quad\dz^\vee=-\dz\,.
$$
Using this anti-involution one can turn a right module $E$ over $\CC[z,\dz]$
into a left module by defining
$$
a\xi\,=\,\xi a^\vee\qquad\textrm{for}\quad \xi\in E\,,\;a\in\CC[z,\dz]\,.
$$
The same applies, of course, to $\CC[t,\dt]$.

\

\pph{\IS-model.}\label{IS-model}
For the \IS-model we take the free module $E$ with basis $\xi_0,\ldots,\xi_n$ over the ring of
Laurent polynomials $\CC[z,z^{-1}]$ and give it the structure of a right module
over $\CC[z,\dz]$ by defining
$$
(\xi_0,\ldots,\xi_n)\dz = (\xi_0,\ldots,\xi_n)(\II-\tilde{A}_1)z^{-1}\,,
$$
where $\II$ is the identity matrix of size $n+1$ and
$\tilde{A}_1=( a_{ij}\,z^{j-i+1})$.

This definition implies
$$
(\xi_0,\ldots,\xi_n)(-z\dz-\tilde{A}_1)\,=\,
(\xi_0,\ldots,\xi_n)(-\dz\,z+(\II-\tilde{A}_1))\,=\,
(0,\ldots,0)\,.
$$
So, Proposition \ref{prop:CH} implies
$$
\xi_0\detr (-z\dz-\tilde{A}_1^\tau)=0\,.
$$
Passing from right to left modules with the involution ${}^\vee$ and
also using Proposition \ref{anti-involution} we can rewrite the above formulas
as
\begin{eqnarray*}
&&\hspace{-2em}\textbf{\IS-model connection:}\quad
(\dz\xi_0,\ldots,\dz\xi_n) = (\xi_0,\ldots,\xi_n)\,(\tilde{A}_1-\II)z^{-1}\,,
\\
&&\hspace{-2em}\textbf{\IS-model differential equation:}\qquad
\detr (\dz z -\tilde{A}_1)\xi_0=0\,.
\end{eqnarray*}

\

\pph{\IS-model (bis).}\label{IS-model(bis)}
We introduce the diagonal matrix
$$
\TT = \mathrm{diag} (0,1,\ldots,n)\,.
$$
Then $z^{\TT}=\mathrm{diag} (1,z\ldots,z^n)$ and
$\tilde{A}_1 = z^{-\TT}\,A\,z^{\TT}z$.
Thus, by the change of coordinates
$(\zeta_0,\ldots,\zeta_n)=(\xi_0,\ldots,\xi_n)z^{-\TT}$ one can
put the above \IS-model connection in the format of \cite{Ma} p. 53 formula (1.23):
$$
(\dz\zeta_0,\ldots,\dz\zeta_n)=(\zeta_0,\ldots,\zeta_n)
(A-(\II+\TT)z^{-1})\,.
$$

\

\pph{\RS-model.}\label{RS-model}
For the \RS-model we take the free module $E'$ with basis $\eta_0,\ldots,\eta_n$ over the ring $\CC[t,\chi_A^{-1}]$, where $\chi_A=\det(A-t)$ is the characteristic polynomial of the matrix $A=(a_{ij})$.
We give $E'$ the structure of a right module
over $\CC[t,\dt]$ by defining
$$
(\eta_0,\ldots,\eta_n)\dt = (\eta_0,\ldots,\eta_n)\,\TT\,(A\,-\,t)^{-1}\,.
$$
Let $\tilde{A}=( a_{ij}\,\dt^{j-i+1})$ and
$V=\mathrm{diag} (1,\dt,\dt^2,\ldots,\dt^n)$.
Then $V\tilde{A}\,=\,\dt A V$ and
$V\dt t\,=\,(\dt t+\TT)V$ and
$$
(\eta_0,\ldots,\eta_n)\,(V\tilde{A}\,-\,V\dt t)=
(\eta_0,\ldots,\eta_n)\,(\dt\,(A\,-\,t)-\TT)V=
(0,\ldots,0)\,.
$$
Thus, if we define the elements $\xi'_0, \ldots, \xi'_n$ in $E'$ by
$\xi'_k\,=\,\eta_k (\dt)^k$, then
$$
(\xi'_0, \ldots, \xi'_n)(\dt t\,-\,\tilde{A})\,=\,(0,\ldots,0)\,.
$$
So, Proposition \ref{prop:CH} implies
$$
\xi'_0\detr (\dt t-\tilde{A}^\tau)=0\,.
$$
Passing from right to left modules with the involution ${}^\vee$ and
also using \ref{anti-involution} we can write this formula also
as
$$
\detr (-t\dt -\tilde{A}^\vee)\eta_0=0\,.
$$
Note $\tilde{A}^\vee\,=\,( a_{ij}\,(-\dt)^{j-i+1})$.
Proposition \ref{diag-change-signs} enables one to
 rewrite the above formulas as
\begin{eqnarray*}
&&\hspace{-1.7em}\textbf{\RS-model connection:}\quad
(\dt\eta_0,\ldots,\dt\eta_n) = (\eta_0,\ldots,\eta_n)\,\TT\,(t\,-\,A)^{-1},
\\
&&\hspace{-1.7em}\textbf{\RS-model differential equation:}\qquad
\detr (t\dt -\tilde{A})\eta_0=0\,.
\end{eqnarray*}
This \RS-model connection has exactly the form of
\cite{Ma} p. 53 formula (1.22).

\

%\pph{Remark.}\label{can be irregular}
%In Section \ref{section:monodromy} we will see that the \RS-model connection
%has only \emph{regular singularities if the matrix $A$ is diagonalizable}.
%On the other hand Example \ref{example irregular} shows that there can be irregular singularities if $A$ is not diagonalizable.
%

\

\pph{Remark.}\label{reducible}
Note that the top row of $\TT$ is zero and that therefore $\eta_0$ is actually
absent from the right hand side of the \RS-model connection formula. Consequently,
the $\CC[t,\chi_A^{-1}]$-submodule $E''$ of $E'$ with basis
$\eta_1,\ldots,\eta_n$ is stable under the action of $\dt$, i.e.
$\dt E''\subset E''$. Moreover, $\dt\eta_0\in E''$ and this implies that
the class of $\eta_0$ modulo $E''$ is horizontal for the connection on
the quotient $E'/E''$. The counterpart of this structure is that the \RS-model differential operator is divisible on the right by $\dt$; see \ref{lemma:DN1}.

\

\pph{Connection vs. differential equation.}\label{horizontal}
We remind the reader how the problem of finding horizontal sections for a
connection relates to that of solving a differential equation.

 Let $\cA$ be a commutative algebra  with a left action of $\CC [t,\dt]$
that satisfies the Leibniz rule. Assume also that $\cA$ acts on the right on
the free $\cA$--module \mbox{$E \otimes_{\CC [t,\chi_A^{-1}]}\cA$} with the Leibniz rule
$$\dt (e a)=(\dt e) a + e (\dt a), \;\; e \in E
\otimes_{\CC [t,\chi_A^{-1}]}\cA,\;
a \in \cA.$$
Let $\Phi$ be an invertible matrix of size $n+1$ whose entries are
elements of $\cA$, and assume that
$$\dt((\eta_0,\ldots,\eta_n)\Phi)\,=\,(0,\ldots,0).$$

(In practice, $\cA$ would of course be the algebra of analytic functions
on some open subset $U$ of $\CC$;
$(\eta_0,\ldots,\eta_n)\Phi$ is then a basis for the space of horizontal analytic
sections
over $U$ for the connection.)

From
\begin{eqnarray*}
&&\dt((\eta_0,\ldots,\eta_n)\Phi)\,=\,(\dt(\eta_0,\ldots,\eta_n))\Phi
+(\eta_0,\ldots,\eta_n)\dt\Phi\\
&&=\,(\eta_0,\ldots,\eta_n)(\TT\,(t\,-\,A)^{-1}\Phi+\dt\Phi)\,=\,
(0,\ldots,0)
\end{eqnarray*}
we see that $\Phi$ is `a fundamental solution matrix' of the system
$$
\dt\Phi\,=\,\TT\,(A\,-\,t)^{-1}\Phi\,.
$$

The trivial fact $\dt(\Phi^{-1}\Phi)=0$ implies
$\dt\Phi^{-1}=\Phi^{-1}\,\TT\,(t\,-\,A)^{-1}$. One can therefore apply the arguments in \ref{RS-model} to the columns of $\Phi^{-1}$ in place of
$(\eta_0,\ldots,\eta_n)$. Thus the left-most column of
$\Phi^{-1}$ is componentwise annihilated by the differential operator
$\detr (t\dt -\tilde{A})$. So the elements appearing as entries in the
left-most column of $\Phi^{-1}$ are solutions of the differential equation
$\detr (t\dt -\tilde{A})f=0$.

\section{Differential operators of type $DN$}
\label{section:DN equations}

\pph{Lemma.}  \label{lemma:DN1} {\it
Let $A=(a_{ij})_{0\lqs i,j\lqs n}$ be a matrix with entries in $\CC$
such that $-A$ is almost triangular. Let $\tilde{A}$ denote the matrix with $(i,j)$-entry $a_{ij}\dt^{j-i+1}$.
Then the differential operator $\detr (t\dt -\tilde{A})$ in the
\RS-model differential equation is uniquely divisible from the right by $\dt$.}

\

\prf
This follows from Proposition \ref{detright2} and the observation that in the matrix $t\dt-\tilde{A}$ the entry in position $(i,j)$ is uniquely divisible from the right by $\dt$ if $i\lqs j$.

\qd

\

\pph{Definition.} \it \label{def:DN1}
With $A$ and $\tilde{A}$ as in Lemma \ref{lemma:DN1} we define the differential operator
$$
L_{A,\infty}\,=\,\detr ( t\dt  -\tilde{A} )\dt^{-1} \,.
$$

\

\pph{Proposition.}  \label{prop:DN1} {\it
The operators $L_{A,\infty}$ one obtains with Definition \ref{def:DN1}
are precisely the operators of the form
$$
(t\dt)^{n}t + \sum_{p = 1}^{n + 1} g_p(t\dt )\,\dt^{p-1}\,.
$$
with $g_p$ a polynomial of degree $\lqs n-p+1$ in $t\dt$.}

\

\prf Setting $u=\dt$ and $\us=t$ in the right-hand side of
\ref{expansion} and applying the anti-involution ${}^\vee$ gives
\begin{eqnarray*}
L_{A,\infty}&=&(-1)^n\left(
t(\dt t)^{n} + \sum_{p = 1}^{n + 1} \sum_{k=0}^{n-p+1}
x^{(p)}_k \dt^{p-1} (\dt t)^k\right)^\vee\\
&=&(t\dt)^{n}t+ \sum_{p = 1}^{n + 1} \sum_{k=0}^{n-p+1}
(-1)^{n+k+p-1}x^{(p)}_k  (t\dt)^k\dt^{p-1}
\end{eqnarray*}
\qd

\

\pph{}\label{DN1 solutions}
Continuing the discussion in \ref{horizontal} we see that the entries in the left-most column of the matrix $\dt \Phi^{-1}$
are solutions to the differential equation $L_{A,\infty}g=0$.

\

\pph{}\label{intro DN0}
Let us now consider the ring of differential operators
$\CC[w,w^{-1},\dw]$ on the torus $\mathrm{Spec }\, \CC[w,w^{-1}]$.
The substitution $t=w^{-1}$ transforms $\dt$ into $-w^2\dw$
and $t\dt$ into $-w\dw$. Thus the differential operator $L_{A,\infty}$,
expanded as in Proposition \ref{prop:DN1}, transforms into
$$
(-w\dw)^{n}w^{-1} + \sum_{p = 1}^{n + 1} g_p(-w\dw )\,(-w^2\dw)^{p-1}\,.
$$

\

\pph{Definition.} \it \label{def:DN0}
With the notations as in \ref{intro DN0} we define the differential operator
$$
L_{A,0}\,=\,(w\dw)^{n} + \sum_{p = 1}^{n + 1}(-1)^n g_p(-w\dw )
\,(-w^2\dw)^{p-1}w\,.
$$
This means, according to \ref{intro DN0}, that the substitution $t=w^{-1}$ transforms $L_{A,\infty}$ into
$(-1)^n\,L_{A,0}\,w^{-1}$.

\

\pph{Proposition.} \label{prop:DN0} {\it
The operators $L_{A,0}$ one obtains with Definition \ref{def:DN0}
are precisely the operators of the form
$$
(w\dw)^n + \sum_{p = 1}^{n + 1}
w^p\:G_p(w\dw)\:\prod_{l=1}^{p-1}(w\dw+l)
$$
where $G_p$ is a polynomial of degree $\lqs n+1-p$ in $w\dw$,
related to the polynomial $g_p$ from Proposition \ref{prop:DN1} by:}
$$
G_p(\mathrm{argument})\,=\,(-1)^{n-p+1}g_p(-\mathrm{argument}-p)\,.
$$

\

\prf
Note $g_p(-w\dw )(w^2\dw)^{p-1}\,w\,=\,w^pg_p(-w\dw-p )\,
\prod_{l=1}^{p-1}(w\dw+l)$.
\qd

\

\pph{}\label{adjoint}
The anti-involution ${}^\vee$ on $\CC[t,\dt]$
defined by $t^\vee=t\,,\;\dt^\vee=-\dt$ maps a differential operator
$L=\sum_{i,j} c_{ij}t^i\dt^j$ to
$
L^\vee\,=\,\sum_{i,j} c_{ij}(-\dt)^j t^i\,.
$
Thus, $L^\vee$ is the \emph{adjoint} of $L$ in the sense of
\cite{Ka} p. 55.
We want to determine the adjoint of the operator
$L_{A,\infty}$.
Recall from Definition \ref{def:DN1} that
$$
L_{A,\infty}\,=\,\detr ( t\dt  -\tilde{A} )\dt^{-1} \,,
$$
where $\tilde{A}$ is the matrix with $(i,j)$-entry
$a_{i j} \dt^{j - i + 1}$. Using Proposition \ref{anti-involution} we find:
$$
L_{A,\infty}^\vee\,=\,-\dt^{-1}\detr (-\dt t-\tilde{A}^{\tau\vee})\,,
$$
where $\tilde{A}^{\tau}$ is the matrix with $(i,j)$-entry
$a_{n-j,n-i} \dt^{j-i + 1}$.

Using Proposition  \ref{diag-change-signs}, one sees
$$
L_{A,\infty}^\vee\,=\,(-1)^n\detr (t\dt -\tilde{A}^\tau)\dt^{-1}
\,=\,(-1)^n
L_{A^\tau,\infty}\,.
$$
Definition \ref{def:DN0} now gives
$$
L_{A,0}^\vee\,=\,(-1)^n wL_{A^\tau,0}\,w^{-1}\,.
$$

\

\pph{Theorem.}\label{self-adjoint} {\it
The following three statements are equivalent
\begin{enumerate}
\item $\quad A=A^\tau$,
\item $\quad L_{A,\infty}^\vee=(-1)^n L_{A,\infty}$,
\item $\quad L_{A,0}^\vee=(-1)^n wL_{A,0}\,w^{-1}$.
\end{enumerate}
}

\

\prf
This follows directly from \ref{adjoint} and \ref{reconstruction}.
\qd

\

\pph{Definition.} \it \label{def:DN--}
Following \cite[2.10]{Go} we call $L_{A,0}$ a
\emph{differential operator of type $DN_{0,0}$} if the matrix $A$ satisfies the condition $A=A^\tau$.
Under the same symmetry condition $L_{A,\infty}$ is called a
\emph{differential operator of type $DN_{\infty,1}$}.
Here $N=n$ is the order of the operator. Thus, for $n=3$ the operator
$L_{A,\infty}$ is of type $D3_{\infty,1}$.

\

\pph{Theorem.}\label{thm:DN--} {\it
\begin{enumerate}
\item
An operator $L$ of type $DN_{\infty,1}$ satisfies $\quad L^\vee=(-1)^nL$.
\item
The operators of type $DN_{\infty,1}$ are precisely the differential operators
of the form
$$
(t\dt)^{n}t + \sum_{p = 1}^{n + 1} g_p(t\dt )\,\dt^{p-1}
$$
with $g_p$ a polynomial of degree $\lqs n-p+1$, that satisfies
$$
g_p(\mathrm{argument})=(-1)^{n-p+1}g_p(-\mathrm{argument}-p).
$$
\item
An operator $L$ of type $DN_{0,0}$ satisfies $\quad L^\vee=(-1)^nwLw^{-1}$.
\item
The operators of type $DN_{0,0}$ are precisely the differential operators
of the form
$$
(w\dw)^n + \sum_{p = 1}^{n + 1}
w^p\:G_p(w\dw)\:\prod_{l=1}^{p-1}(w\dw+l)
$$
where $G_p$ is a polynomial of degree $\lqs n-p+1$, that satisfies
$$
G_p(\mathrm{argument})=(-1)^{n-p+1} \:G_p(-\mathrm{argument}-p)\,.
$$
\end{enumerate}
}

\

\prf
(i) follows from Theorem \ref{self-adjoint} and Definition
\ref{def:DN--}.
\\
(ii) Write the operator $L$ as in Proposition \ref{prop:DN1}.
Then $L^\vee=(-1)^nL$ becomes
$$
t(-\dt t)^{n} + \sum_{p = 1}^{n + 1} (-\dt)^{p-1}g_p(-\dt t )
\,=\,
(-1)^n\,\left[(t\dt)^{n}t + \sum_{p = 1}^{n + 1} g_p(t\dt )\,\dt^{p-1}\right]
$$
and boils down to
$$
g_p(t\dt)\dt^{p-1}\,=\,(-1)^{n+p-1}\dt^{p-1}\,g_p(-\dt t)\,=\,
(-1)^{n+p-1}
g_p(-t\dt -p)\,\dt^{p-1}
$$
for every $p\gqs 1$. This proves (ii).

(iii) and (iv)
follow directly from Theorem \ref{self-adjoint}, (ii) and Proposition \ref{prop:DN0}.
\qd

\

\pph{Example.}\label{example irregular}
For $\gl\in\CC$ consider the matrix
$$
A=\left(\begin{array}{rrr}\gl&-\frac{3}{2}\gl^2&-\gl^3\\
1&-2\gl&-\frac{3}{2}\gl^2\\ 0&1&\gl\end{array}\right)
$$
A straightforward calculation shows that for this matrix
$$
L_{A,\infty}\,=\,t^3\dt^2\,+\,3t^2\dt\,+\,t-\gl\,.
$$

Divided by $t$, it yields a monic polynomial in
$t\dt$ with coefficients in $\CC(t).$ In the case $\gl\neq 0$  one of these coefficients,
$\lambda/t,$ is not analytic at 0   hence, by the Fuchs criterion,
 the operator $L_{A,\infty}$ has an
irregular singularity at $t=0$.

\

\pph{Remark.}\label{remark irregular}
The above example shows that operators of type $DN_{\infty,1}$ may
have irregular singularities.
On the other hand we will see in the next section
that if the matrix $A$ is diagonalizable the corresponding
operators $DN_{\infty,1}$ and $DN_{0,0}$ have
only regular singularities.

\section{Monodromy}
\label{section:monodromy}

\pph{Assumption.}\label{assumptions}
\emph{In this section we assume that
$-A$ is almost triangular, $A$ is diagonalizable and $A=A^\tau$.}

\

\pph{}\label{fuchs system}
We are going to investigate the local monodromies of the \RS-model connection,
$$
(\dt\eta_0,\ldots,\dt\eta_n) = (\eta_0,\ldots,\eta_n)\,\TT\,(t\,-\,A)^{-1},
$$
or equivalently (cf. \ref{horizontal}) of the system
$$
\dt\Phi\,=\,\TT\,(A\,-\,t)^{-1}\Phi\,,
$$

It is clear that the singularities are at $\infty$ and at the eigenvalues
of $A$.
Since $-A$ is almost triangular, the last coordinate
of every eigenvector of $A$ is non-zero (if it were zero, then the value of the
pairing of the bottom row  of $A$ with the eigenvector would be zero, hence
so would be the $(n-1)$th coordinate, and so forth). This implies that all
eigenspaces of $A$ have dimension $1$. Since $A$ is diagonalizable, this means
that all eigenvalues of $A$ have multiplicity $1$.
Let
$$
A= C\gL C^{-1}\,,\qquad\gL=\mathrm{diag}(\gl_0,\ldots,\gl_n)\,.
$$
By assumption \ref{assumptions}  and by \ref{def:transpose}
we have $A=A^\tau=\JJ A^t\JJ$. Hence,
$C\gL C^{-1}=(C^t\JJ)^{-1}\gL C^t\JJ$.
Since $\gl_i\neq\gl_j$ if $i\neq j$, this
implies that $C^t\JJ C$ is a diagonal matrix.
By multiplying $C$ on the right by a suitable diagonal
matrix, we may assume
$$
C^t\JJ C= \II\,.
$$
Let $\bu_0,\ldots,\bu_n$ be the columns of the matrix $C$.
So, $\bu_i$ is an eigenvector of $A$ for the eigenvalue $\gl_i$.

\

\pph{Theorem.}\label{fuchs} {\it
With the assumptions and notations of \ref{assumptions} and \ref{fuchs system},
\begin{enumerate}
\item the \RS-model system can be written as
$$
\dt\Phi\,=\,\sum_{j=0}^n\frac{1}{t-\gl_j}S_j\Phi
$$
where $S_j=-\TT C E_j C^{-1}$ and $E_j$ is the $(n+1)\times(n+1)$-matrix with
$1$ in position $(j,j)$ and zeros elsewhere;

\item the vectors $\bu_i$ with $i\neq j$ are eigenvectors of $S_j$
for the eigenvalue $0$ and $\TT\bu_j$ is an eigenvector of $S_j$
for the eigenvalue $-\frac{n}{2}$.
\end{enumerate}}

\

\prf
The first statement is an immediate consequence of
$$
(A-t)^{-1}\,=\,C(\gL-t)^{-1}C^{-1}=\sum_{j=0}^n\frac{1}{\gl_j-t}C E_j C^{-1}\,.
$$
From $S_jC=-\TT C E_j$ one sees that
the vectors $\bu_i$ with $i\neq j$ are eigenvectors of $S_j$
for the eigenvalue $0$ and that $\TT\bu_j$ generates the image of $S_j$.
So $\TT\bu_j$ is also an eigenvector of $S_j$ with eigenvalue equal to $\mathrm{trace}(S_j)\,=\,\mathrm{trace}(-C^{-1}\TT C E_j)$. The remaining eigenvalue thus equals
the $(j,j)$ entry of the matrix $-C^{-1}\,\TT\,C$.
Since $(C^{-1}\,\TT\,C)^t=C^t\,\TT\,(C^t)^{-1}=C^{-1}\,\JJ\TT\JJ\,C$,
the diagonal of $C^{-1}\,\TT\,C$ equals $\frac{1}{2}$ times
the diagonal of $C^{-1}\,(\TT+\JJ\TT\JJ)\,C\,=\,n\II$.
Therefore $\mathrm{trace}(S_j)=-\frac{n}{2}$.
\qd

\

\pph{Theorem.}\label{monodromies} {\it
Under the assumption \ref{assumptions} the singularities of the system of differential equations
$$
\dt\Phi\,=\,\TT\,(A\,-\,t)^{-1}\Phi
$$
are regular singular points located at $\infty$ and at the (distinct) eigenvalues
 $\gl_0,\ldots,\gl_n$ of $A$.

The monodromy transformation $M_\infty$ along a small positively oriented
simple loop around $\infty$ is maximally unipotent, i.e. $(M_\infty-\II)^n\neq 0$,
$(M_\infty-\II)^{n+1}=0$.

In case $n$ is odd
the monodromy transformation $M_j$ along a small positively oriented simple
loop around $\gl_j$ has an eigenvalue
$1$ with $n$-dimensional eigenspace and
an eigenvalue $-1$ with $1$-dimensional eigenspace.

In case $n$ is even $1$ is the only eigenvalue of
the monodromy transformation $M_j$ along a small positively oriented
simple loop around $\gl_j$ and the dimension of the eigenspace is $\gqs n$.
}

\

\prf The proof of this theorem is given in \ref{fuchsianities},
\ref{infinity monodromy} and \ref{finite monodromy}.

\

\pph{}\label{fuchsianities} \bf Formal Fuchsian theory. \rm
We borrowed the following concise account of it from
 \cite[Ch.III \S 8]{DGS}.
Let $G$ be a square matrix with entries in $\CC[[x]]$.
Write $\dx=x\frac{d}{dx}$.
 For
an invertible matrix $H$ with entries in $\CC[x,x^{-1}]$, define
$G_{[H]}$ as $G_{[H]}=(\dx H) H^{-1}+HGH^{-1}.$
If a matrix
$\Phi$ of functions in $x$ satisfies $\dx\Phi=G\Phi$, then
$\dx (H\Phi)=G_{[H]}(H\Phi)$.

Let $\alpha$ be an eigenvalue of $G(0)$. Then according to \cite[Lemma
III.8.2]{DGS} there is an invertible matrix $H$ with entries in $\CC[x,x^{-1}]$
such that the matrix $G_{[H]}$ has entries in $\CC[[x]]$ and
such that the matrix $G_{[H]}(0)$ has the same eigenvalues, counted with
multiplicities, as $G(0)$, except that
$\alpha$ is replaced by $\alpha+1$. A similar construction with $H^{-1}$
instead of
$H$ replaces $\alpha$ by $\alpha-1$.
By repeated application of this construction we find an invertible matrix $H$ with entries in $\CC[x,x^{-1}]$
such that the matrix $G_{[H]}(0)$ has \emph{prepared eigenvalues}; i.e.
the eigenvalues $\alpha_0,\ldots,\alpha_n$ satisfy: if $\alpha_i\in\ZZ$ then $\alpha_i=0$,
and if $\alpha_i-\alpha_j\in\ZZ$ then $\alpha_i=\alpha_j$.
The system $\dx\Psi=G_{[H]}\Psi$ then has a solution matrix of the form
$\Psi=W\exp(G_{[H]}(0)\log x)$ where $W$ is a matrix with entries in $\CC[[x]]$
such that $W(0)$ is the identity matrix.
Thus if $x$ runs through a small positively oriented loop about $0$ the fundamental solution matrix $\Phi=H^{-1}\Psi$ gets multiplied from the right by $\exp(2\pi i G_{[H]}(0))$.

\

\pph{}\label{infinity monodromy}
For the local monodromy of our system
$
\dt\Phi\,=\,\TT\,(A\,-\,t)^{-1}\Phi
$
around $\infty$ we use the local coordinate
$x=t^{-1}$. This puts the system
 in the form
$$
x\frac{d}{dx}\Phi\,=\,G\Phi \qquad\textrm{with}\qquad G=\TT\,(\II-Ax)^{-1}\,.
$$
So, $G(0)=\TT=\mathrm{diag} (0,1,\ldots,n)$. Take
$H=\mathrm{diag} (1,x^{-1},x^{-2},\ldots,x^{-n})$.
Then
$$
G_{[H]}=(\dx H) H^{-1}+HGH^{-1}=\TT(-\II+\sum_{k\gqs 0}(xH\,A\,H^{-1})^k)\,.
$$
We see that $G_{[H]}(0)$ is a nilpotent matrix with
$G_{[H]}(0)^n\neq 0$ and $G_{[H]}(0)^{n+1}=0$.
So, $G_{[H]}(0)$ has prepared eigenvalues and the monodromy is given by
the maximally unipotent matrix $\exp(2\pi i G_{[H]}(0))$.

\

\pph{}\label{finite monodromy}
For the local monodromy around an eigenvalue $\lambda_j$ of $A$
we apply the general theory in \ref{fuchsianities} with $x=t-\lambda_j$ and $G=\TT\,x(A\,-\,\lambda_j\II-x\II)^{-1}$.
Then $G(0)=S_j$ as in Lemma \ref{fuchs}.

Thus, if $n$ is odd $G(0)$ has prepared eigenvalues. The matrix
for the monodromy around $\gl_j$ is $\exp(2\pi i G(0))$  and has an eigenvalue
$1$ with $n$-dimensional eigenspace and
an eigenvalue $-1$ with $1$-dimensional eigenspace.

If $n$ is even $G(0)$ does not have prepared eigenvalues. The general theory now provides a matrix $H$ such that $0$ is the only eigenvalue of
$G_{[H]}(0)$ and the eigenspace has dimension $\gqs n$.

\text{This completes the proof of Theorem \ref{monodromies}}.

\

\pph{}\label{DN monodromies}
Let us now turn to the local monodromies of the operator $L_{A,\infty}$
defined in \ref{def:DN1}. As we have seen in \ref{horizontal}
the entries in the left-most column of the matrix $\Phi^{-1}$
(where $\Phi$ is as in Theorem \ref{monodromies}) constitute a basis for
the solution space of the differential equation
$\detr (t\dt -\tilde{A})f=0$.
So the local monodromies act on this basis by multiplication from the left
by certain matrices $M_\infty^{-1}$ and $M_j^{-1}$ ($j=0,\ldots,n$).

We know from \ref{def:DN1} that
$\detr (t\dt -\tilde{A})\,=\,L_{A,\infty}\dt$.
So the constant functions are solutions of the differential equation
$\detr (t\dt -\tilde{A})f=0$ and the solution space of the differential equation $L_{A,\infty}g=0$ is obtained by taking the solution space of
$\detr (t\dt -\tilde{A})f=0$ modulo the constant functions.
Since the constant functions are fixed by all local monodromies we derive from
Theorem \ref{monodromies} the following result about the monodromies of the
operator $L_{A,\infty}$:

\

\pph{Corollary. }\label{cor:monodromies} {\it
Under the assumption \ref{assumptions} the differential equation
$$
L_{A,\infty}g=0
$$
has regular singularities located at $\infty$ and at the (distinct) eigenvalues
 $\gl_0,\ldots,\gl_n$ of $A$.

The local monodromy around $\infty$ is maximally unipotent,
i.e. represented by an
$n\times n$-matrix $\tilde{M}_\infty$ satisfying
$(\tilde{M}_\infty-\II)^{n-1}\neq 0$, $(\tilde{M}_\infty-\II)^n=0$.

If $n$ is odd
the local monodromy around $\gl_j$ has an eigenvalues
$1$ and $-1$ with eigenspaces of dimension $n-1$ and $1$, respectively.

If $n$ is even $1$ is the only eigenvalue of
the local monodromy around $\gl_j$ and the dimension of the eigenspace is
$\gqs n-1$.}
\qd
%

%%%

\section{Polarizability and generic irreducibility}

In this section we prove that a generic operator of type $DN$ is
irreducible. The tactic is as follows. First we prove that
reducibility is a closed condition on an open set of the affine
space of parameters of $DN$s. Up to this point we deal with the
$DN_{\infty,1}$ flavor. Then we exhibit a single irreducible $DN$
operator $\mathcal{H}$, for aesthetic reasons in the $DN_{0,0}$
form, such that any operator in its neighborhood in the analytic
topology is still irreducible.

We start, however, with the polarizability theorem:

\

\pph{Theorem}. \it The monodromy of a differential operator $L$ of
type DN is polarized (i. e. its monodromy representation respects a
non-degenerate bilinear form).

\

\prf  The operator $L$ determines a vector bundle over
$\mathbb{A} ^1 (\CC)  \setminus \{\text{singularities of } L \}$ and
a connection in it. The adjoint gives rise to the dual connection in
the dual bundle. Hence, the monodromy representation that
corresponds to the adjoint is contragredient to the original one. We
have proved in Theorem \ref{self-adjoint} that the differential
operator $L$ of type $DN_{\infty,1}$  coincides, up to a sign, with
its adjoint. Hence, its monodromy representation is isomorphic to
its contragredient, and any such isomorphism corresponds to a
non-degenerate bilinear form respected by the monodromy
representation.
\qd

\

If a polarized monodromy representation is irreducible, then
the polarizing bilinear form is defined uniquely up to a scalar
(being an endomorphism of an irreducible object), and is either
symmetric or skew depending on how the argument interchange
involution acts on its span.

\bigskip

In our discussion of irreducibility we consider differential operators
with rational coefficients, i.e.
elements of $\CC(t)[\partial_t]$. As before, $N=n$.

\bigskip

\pph{Definition.} \it  Let $L=c_n \partial ^n + c_{n-1} \partial
^{n-1}+ \dots + c_0 \in \CC (t)[\partial],$ and $p$ be a point in
$\PP^1(\CC)$. The residue $\res _p L$ is by definition $\res_p
(c_{n-1}/c_n) dt.$

\
%
%\pph{Proposition.} The sum of the residues of $L$ over $\PP ^1(\CC)$
%is zero.

\

\pph{Proposition.} \label{residues} {\it  Let $L=c_n \partial ^n +
c_{n-1} \partial ^{n-1}+ \dots + c_0 $  be a differential operator
of type $DN_{\infty,1}$ with $n+1$ distinct finite singularities.
Then:
\begin{enumerate}
\item  The residue of
$L$ at any finite singularity $p$ is $\frac{n}{2}.$
\item  The residue of
$L$ at infinity  is $-n(n+1)/2$.
\end{enumerate}
}

\

\prf Recall that we have identified $L$'s of type $DN_{\infty,1}$ in
Theorem \ref{thm:DN--} ii) as the differential operators of the form

$$
(t\dt)^{n}t + \sum_{p = 1}^{n + 1} g_p(t\dt )\,\dt^{p-1}
$$
with $g_p$ a polynomial of degree $\lqs n-p+1$, that satisfies
$$
g_p(\mathrm{argument})=(-1)^{n-p+1}g_p(-\mathrm{argument}-p).
$$

%$$
%(w\dw)^n + \sum_{p = 1}^{n + 1}
%w^p\:G_p(w\dw)\:\prod_{l=1}^{p-1}(w\dw+l)
%$$
%where $G_p$ is a polynomial of degree $\lqs n-p+1$ that satisfies
%$G_p(\mathrm{argument})=(-1)^{n-p+1} \:G_p(-\mathrm{argument}-p)\,$.
%

Consider the $p$-th term to the right of the summation sign in the
formula above. It can be presented as
$$ \bar g_p(t\dt+\frac{p}{2}){\dt}^{p-1},$$
where $\bar g_p$ is, depending on the parity of $n-p+1$, an even or
odd polynomial (i.e. as a function of its argument). Let  $\sum
c^{(p)}_i(t) \dt^i$ be the expansion of the $p$-th term.

To compare  $c^{(p)}_{n}$ with  $c^{(p)}_{n-1}$ it suffices to
assume that  $\bar g_p$ is a monomial (as a function of its
argument) of the degree $n-p+1$: the lower order terms contribute
neither to $c^{(p)}_{n}$ nor to $c^{(p)}_{n-1}$. Inductively on $l$,
$$(t\dt+\frac{p}{2})^l=t^l \dt^l+\frac{l(l+p-1)}{2}t^{l-1}\dt^{l-1}+ \dots ,$$ and we compute
$c^{(p)}_{n}=t^{n-p+1}$ and
$\displaystyle{c^{(p)}_{n-1}=\frac{(n-p+1)n}{2}t^{n-p}}$,
hence
$c^{(p)}_{n-1}= \displaystyle{\frac{n}{2} \frac{d c^{(p)}_{n}}{dt}}$.
Arguing in the same vein about
$(t\dt)^{n}t=t(t\dt+1)^{n}$, we see, by linearity, that $c_{n-1}=
\displaystyle \frac{n}{2}{\frac{dc_n}{dt}}.$
Expanding $c_n$ into a series at a singularity proves the assertion (i).

Assertion (ii) follows because the sum of the residues of a rational function is $0$.

\qd

\

\pph{Definition.} \it  We say that a differential operator $L=c_n
\partial ^n + c_{n-1} \partial ^{n-1}+ \dots + c_0 \in \CC
(t)[\partial]$ is irreducible if it cannot be represented as a
product $L=L_1L_2$ with $L_1,L_2 \in \CC (t)[\partial]$ of positive
order.

\

 \rm Consider the affine space $\mathbb{A}=\Spec \, \CC [a_{ij}], \;   0 \lqs i \lqs
j \lqs n  , $ of differential operators of type DN.

\

 \pph{Proposition}. \label{ir-locally-closed} {\it The locus of reducible
 DN operators is closed in a non-empty
open subset of $\mathbb{A}$.

\

\prf Let $\{ L= c_n\dt^n + c_{n-1} \dt^{n-1} + \dots + c_0 \}$ be
any set of differential operators  with rational coefficients such
that:
\begin{enumerate}
\item the degrees of numerators and denominators of all $c_i$ are bounded.
\item there is a finite set $R \subset \CC$ such that the residue at
any singularity of any operator in the family is in $R$.
\end{enumerate}
According to \cite[section 9]{vH}, there exists a positive integer
$h$ such that for any factorization
$$L=L'L''  \eqno \text{(\bf F\rm)} $$
where $L'= \dt^{n'} + c'_{n'-1} \dt^{n'-1} + \dots + c_0' $ and
$L''= c''_{n''}\dt^{n''} + c''_{n''-1} \dt^{n''-1} + \dots + c_0'' $
with $c_i',c_i'' \in \CC(t),$ the degrees of numerators and
denominators of all $c_i', c_i''$ are bounded from above by $h.$

According to Proposition \ref{residues}, the set $D.S.$ of all
$DN$'s with $n+1$ distinct finite singularities satisfies hypotheses
i) and ii) above. Therefore, the factorization (\bf F\rm) can be
interpreted as defining a subscheme of the product of (an open
subscheme of) the affine space $\mathbb{A}$ of the parameters
$a_{ij}$ of $DN$ and (open subschemes of) affine spaces of variable
coefficients of the numerators and denominators of $c'$'s and
$c''$'s. By Chevalley's theorem\footnote{Under a morphism of
finite type of Noetherian schemes, the image of a finite  union of
locally closed subsets is a finite  union of locally closed subsets.
- \cite[Ex. 3.19]{Ha}.} the reducible locus is a finite union of
locally closed subsets of $D.S.$ Hence, reducibility is a closed
condition on an open subset of $\mathbb{A}$.
\qd

\

\pph{} Consider the operator
$$\mathcal{H}=D^n-w^{n+1}(D+1)(D+2)\dots (D+n).$$
Up to a scalar, it is the pullback of the regular hypergeometric
operator
$$\mathcal{H}'=D_u^n-u(D_u+\frac{1}{n+1})(D_u+\frac{2}{n+1})\dots
(D_u+\frac{n}{n+1})
$$
under the Kummer cover $u=w^{n+1}.$ The coefficients of
$\mathcal{H}'$ at $u^0$ and $u^1$ are the polynomials $D_u^n$ and
$(D_u+\frac{1}{n+1})(D_u+\frac{2}{n+1})\dots (D_u+\frac{n}{n+1})$
that have no common roots $ \text{mod } \ZZ$.
 According to the classical criterion \cite[3.2.1]{Ka} its monodromy is irreducible.
 It can be computed as in \cite{BH}, or using the following proposition
 from \cite{RRV}:

 \

\pph{Proposition} \cite[1.2.2]{RRV}. Let $F$ be a linear space
endowed with a non-degenerate symmetric (or skew) form $(\,, \,)$.
Let $U$ be a non-degenerate operator acting on $F$. Let $v$ be a
cyclic vector for the operator $U$ and let $S$ be the reflection
with respect to $v$:
$$
S:x \mapsto x-(x,v)v.
$$

Then
$$
1+\sum_{i=1}^\infty(U^iv,v)t^i=\frac{\det(1-tUS)}{\det(1-tU)}.
$$

\qd

\

\pph{} \label{bracket} Consider the global monodromy of
$\mathcal{H}'$. Let $U$ be its monodromy around  $\infty$ and $S$ be
the monodromy around $1$. The monodromy around $0$ is then the
inverse of $US$. Computing for instance by the standard Fuchsian
procedure described in \ref{fuchsianities}, one arrives at a
classically known description of polarized hypergeometric monodromy
\cite[3.2-3.3, 3.4, 3.5.4, 3.5.8]{Ka}. There is a unique, up to a
scalar, non-degenerate bilinear form respected by $U$ and $S$. It is
symmetric for $n$ odd and skew for $n$ even.   $S$ is a
pseudoreflection\footnote{I.e., $S-\II$ has rank 1.} in $O$ (resp.
$Sp$), hence a symmetric (resp. skew) reflection with respect to
some vector $v$. $US$ is unipotent. The eigenvalues of $U$ are all
non-trivial $n+1$-th roots of unity, each with multiplicity $1$.

%In particular, $v$ is a cyclic vector for $U$ (the $\CC[U]$--span
%of $v$ is stable under $U$ and $S$) and the previous
%proposition applies.
%% $\displaystyle{\{ \exp (\frac{2 \pi i j}{n+1}) \mid j= 1,\dots, n+1 \}}$.
%
%The global monodromy of $\mathcal{H}$ is the index $n+1$ subgroup of
%that of $\mathcal{H}'$, generated, say, by the reflections with
%respect to $v, Uv, U^2v, \dots, U^nv$. Computing
%$$(U^jv,
%U^iv)=(U^{j-i}v,v)=\pm \, C_{n+1}^{j-i} \text{ for } 0\le i \ne j
%\le n \eqno (\bf B)$$ as in the proposition above shows that no two
%of $U^iv$ are orthogonal. Hence, the subgroup generated by these
%reflections acts on the fiber $F$ irreducibly.

In particular, $v$ is a cyclic vector for $U$, because the
$\CC[U]$-span
 of $v$ is stable under $U$ and $S$ and because
 the monodromy representation of $\mathcal{H}'$ is irreducible.

 The global monodromy of $\mathcal{H}$ is the index $n+1$ subgroup of
 that of $\mathcal{H}'$, generated, say, by the reflections with
 respect to $v, Uv, U^2v, \dots, U^nv$. Compute  $(U^j v,U^i)$  by
 applying
 Proposition 6.7. :

 $$(U^jv,U^iv)=(U^{j-i}v,v)=\pm \, C_{n+1}^{j-i} \text{ for } 0\le i \ne j
 \le n \eqno (\bf B)$$.

 This shows that no two
 of $U^iv$ are orthogonal. Hence, the subgroup generated by these
 reflections acts on the fiber $F$ irreducibly.

\

\pph{Proposition.} {\it The operator $\mathcal{H}$ is of type $DN_{0,0}$. It
has $n+1$ distinct singularities and is irreducible.}

\

\prf The first assertion follows from the presentation in
\ref{thm:DN--} iv), the second one is obvious. The third one was
proved in the previous paragraph: in a non--trivial factorization
$\mathcal{H}=\mathcal{H}_1\mathcal{H}_2$ both factors need to be
regular singular, in particular, the local system of solutions of
$\mathcal{H}_2 \Phi =0$ would be a non--trivial subsystem of the one
for $\mathcal{H}. $ \qd

\

\pph{Theorem.} \it A generic differential operator $L$  of type DN
is irreducible.

\

\prf We will spell out a semicontinuity
argument in a sufficiently small open neighborhood of $\mathcal{H}$
in analytic topology of $\mathbb{A} (\CC)$ to show that the
monodromy stays irreducible. By \ref{cor:monodromies}, the monodromy
is generated by $n+1$ pseudoreflections of the form $x \mapsto
x+h_i(x)v_i, \; i = 0, \dots, n$ where each $h_i$ stands for a
non-zero covector in a fiber vector space, and each $v_i$, for a
non-zero vector. Let $H_i= \ker h_i.$

Assume the monodromy representation is reducible, and $F_0$ is a
monodromy-stable proper subspace of $F$. An alternative is now
associated with each $v_i$  : \sl $F_0$ either contains $v_i$, or is
in $H_i$. \rm Reduce, if necessary, our neighborhood and note that
$F_0$ cannot contain all $v$'s (resp., be contained in all $H$'s)
because the subspace spanned by all $v$'s is $F$ (resp., the
intersection of all $H$'s is zero), since this is the case for
$\mathcal{H}$. Hence, there would exist $i$ and $j$ such that $v_j
\in H_i$, which again is impossible because it does not happen for
$\mathcal{H}$: no bracket in formula \bf B \rm of \ref{bracket} is
zero.

Summing up,  the monodromy stays irreducible under deformation of
$\mathcal{H}$ in the class of non-degenerate DNs; it implies that
irreducibility is Zariski dense; we showed that reducibility is
locally closed in \ref{ir-locally-closed}. Therefore, the locus of
irreducible DNs contains a non-empty Zariski open set.
\qd

\bigskip

A corollary is

\

\pph{Theorem.} The monodromy representation of a generic
(non-degenerate) operator DN is polarizable by a bilinear form which
is unique up to a scalar. The polarization is skew if $N$ is even
and symmetric if $N$ is odd.

\qd

\bigskip
\bigskip
\bigskip

The authors are grateful to Victor Przyjalkowski for proofreading
the note.

%%%%%%%%%%%%%%%%%%%

%%%%%%%%%%%%%

\end{document}